\documentclass[12pt]{article}
\usepackage{latexsym}
\usepackage{amssymb}

\newcommand{\n}{\mathfrak n}
\newcommand{\T}{\mathfrak T}
\renewcommand{\ll}{\lesssim}
 
\newenvironment{procl}[1]{\medskip\par\noindent {\bf #1.}  \it}%
{\par\noindent}

\title{Variations on the theme of \\  
Marcinkiewicz' inequality}
\author{V. Matsaev
\thanks
{Supported in part by the Israel Science Foundation of the Israel Academy
of Sciences and Humanities under Grant No. 93/97-1.}, I. Ostrovskii
\thanks{Supported in part by the INTAS Project No. 96-0858.},
and M. Sodin$^{\ast}$}

\begin{document}

\maketitle

\section{Introduction}
\setcounter{equation}{0}

Let $h$ be a smooth function on $\bf R$ with a compact support, and 
$$
g(x) = \hbox{p.v.} \frac{1}{\pi}\, \int_{\bf R} \frac{h(t)}{t-x}\, dt
$$
be its Hilbert transform. Set $f=g+ih$ and introduce the function 
$$
u_f(z)=\int_{\bf R} H(zf(t))\, dt\,, \qquad H(z)=\log|1-z|+\hbox{Re}(z)\,, 
$$
called the {\em logarithmic determinant} of genus one.
It is subharmonic in $\bf C$, and its Riesz measure is 
$d\mu_f(\zeta)=d\nu_f(\zeta^{-1})$, where $d\nu_f$ 
is the distribution measure of $f$:
$$
\nu_f(E)=\hbox{meas}\left(\{t:\, f(t)\in E\}
\right)\,, \qquad E \ \hbox{is a borelian subset
of} \ \bf C\,,
$$
and $\hbox{meas}(\,.\,)$ stands for the Lebesgue measure on $\bf R$. 
Let 
$$
\mu_f (r)=\mu_f(\{|z|\le r\})=\hbox{meas}\left(\{|f|\ge r^{-1}\} \right)
$$ 
be a (conventional) counting function of $d\mu_f$,
and let
\begin{eqnarray*}
\n_f(r) &=& 
\mu_f \left( \{ |z - ir/2| \le r/2\}\right) + 
\mu_f \left( \{ |z + ir/2| \le r/2\}\right) \\ 
&=& \mu_f\left(\{ |\hbox{Im}(z^{-1})|\ge r^{-1}\}\right)
= \hbox{meas}\left( \{|h|\ge r^{-1}\} \right)
\end{eqnarray*}
be its {\em Levin-Tsuji counting function}, see 
\cite{Levin1}, \cite{Tsuji}, and \cite[Chapter~1]{GO}. Then the
classical estimates of the Hilbert transform can be easily rewritten as
upper bounds of $\mu_f(r)$ by $\n_f(r)$.   

For example, Marcinkiewicz' inequality (see \cite[Chapter~V]{Koosis})
\begin{equation}
m_f(\lambda) \le 
C \,\left\{
\frac{1}{\lambda^2}\, \int_0^\lambda sm_h(s) ds 
+ \frac{1}{\lambda}\, \int_\lambda^\infty m_h(s) ds 
\right\}\,, \qquad 0<\lambda<\infty\,, 
\label{0.1}
\end{equation}
where 
$$
m_f(\lambda) = \hbox{meas}\left(\{|f|\ge \lambda\}\right)
=\nu_f(\{|w|\ge \lambda\}) = \mu_f (\lambda^{-1})\,, 
$$
and
$$ 
m_h(\lambda) =  \hbox{meas}\left(\{|h|\ge \lambda\}\right)
=\nu_f(\{|\hbox{Im} w|\ge \lambda\}) = \n_f (\lambda^{-1})\,, 
$$
reads:
\begin{equation}
\mu_f(r) \le C
\left\{
r \int_0^r \frac{\n_f(t)}{t^2}\, dt + r^2 \int_r^\infty
\frac{\n_f(t)}{t^3}\, dt
\right\}\,,
\qquad 0<r<\infty\,.
\label{0.2}
\end{equation}
From this, one readily obtains
\begin{equation}
\mu_f(r) \le Cr \int_0^\infty \frac{\n_f(t)}{t^2}\, dt\,, \qquad
0<r<\infty\,,
\label{0.3}   
\end{equation}
which is equivalent to Kolmogorov's weak $L^1$-type inequality:
$\lambda m_f(\lambda) \le C||h||_{L^1}$, $0<\lambda<\infty$,
and
\begin{equation}
\int_0^\infty \frac{\mu_f(t)}{t^{p+1}} dt
\le C(p) \int_0^\infty \frac{\n_f(t)}{t^{p+1}} dt\,,
\qquad 1<p<2\,,
\label{0.5}
\end{equation}
which is equivalent to M.~Riesz' inequality: 
\begin{equation}
||f||_{L^p} \le C(p) ||h||_{L^p}\,.
\label{0.6}
\end{equation}
The classical proof of inequality (\ref{0.1}) is based on the
interpolation technique which later served as one of the cornerstones of
an abstract theory of interpolation of operators in Banach spaces
\cite{Bennett}.

A natural question arises: {\em what is special about the subharmonic
function $u_f(z)$ which makes inequality (\ref{0.2})
valid?} The answer is proposed in \cite{MS1}: a key fact
ensuring this, is the positivity condition:
\begin{equation}
u_f(z) \ge 0\,, \qquad z\in {\bf C}\,,
\label{0.7}
\end{equation}
which can be easily checked with the help of Green's
formula (see \cite[Lemma~5]{Essen}) or by applying of the Cauchy
residue theorem (see \cite{MS2}). 

This leads to a heuristic principle which suggests that
\begin{itemize}
\item[$\bullet$]
{\em results about the distribution of the Hilbert transform can be
deduced from inequality (\ref{0.7}) by using methods of the subharmonic
function theory}.
\end{itemize}
As will be shown, inequality (\ref{0.7}) is even too strong, and in many
cases it suffices to assume that $u_f(x)\ge 0$ on $\bf R$, or to
control the negative part $u_f^-=\max(0,-u_f)$ on $\bf R$.

The principle shows a path to new results on the Hilbert transform. In
\cite{MS2}, its application leads to a complete description of the
distribution of the Hilbert transform of $L^1$-functions and measures of
finite variations. At the same time, putting known estimates of the
Hilbert transform into this setting, we arrive at new interesting
questions about the argument-distribution of the Riesz measure
in certain classes of subharmonic functions.  For example, the proofs 
of the inequalities of Kolmogorov and M.~Riesz 
found in \cite{MS1} give new bounds for zeros of polynomials.
Positivity condition (\ref{0.7}) links our work with the theory of
uniform algebras and Jensen measures (see \cite{Gamelin}).

In this work we put forward a new approach to the Marcin\-kie\-wicz
inequality (\ref{0.1}) (or (\ref{0.2})).
The methods applied in \cite{MS1}, \cite{MS3} are too rigid for this. 
Here we use a different technique.

\medskip
Here and later on, we use the following notations: 

\smallskip\par\noindent $\phi (s) \ll \psi (s)$ means that there is a
positive numerical constant $C$ such that, for each $s>0$, $\phi (s) \le C
\psi (s)$;

\smallskip\par\noindent 
$\phi (s) \ll_{\alpha} \psi (s)$ means the same as above but $C$ may
depend on
a parameter $\alpha$;

\smallskip\par\noindent $H(z)=\log|1-z| + \hbox{Re}(z)$ 
is the canonical kernel of genus one;

\smallskip\par\noindent $\bf C_\pm$ are the upper and lower half-planes.

\section{Main results}
\setcounter{equation}{0}

Define a subharmonic canonical integral of genus one:
\begin{equation}
u(z)=\int_{\bf C} H(z/\zeta)\, d\mu(\zeta)\,, 
\label{1.1}
\end{equation}
where $d\mu$ is a non-negative locally finite measure on $\bf C$ such that
\begin{equation}
\int_{\bf C}
\min\left(\frac{1}{|\zeta|},\frac{1}{|\zeta|^2}\right)\,d\mu(\zeta)<\infty\,.
\label{1.2}
\end{equation}
Let
$$
\mu (r) = \mu (\{|z|\le r\})
$$
be a (conventional) counting function of the measure $d\mu$, and let 
$$
\n(r) = \mu\left( \{ |\mbox{Im}(1/z)|\ge 1/r \} \right) 
= \mu\left( \{|z - ir/2| \le r/2 \} \right)
+ \mu\left( \{|z + ir/2| \le r/2 \} \right)
$$
be its {\em Levin-Tsuji counting function} \cite{Levin1}, \cite{Tsuji} 
(see also \cite[Chapter~1]{GO}).

Let $M(r,u)=\max_{|z|\le r} u(z)$. Then by the Jensen inequality, 
$\mu (r) \ll M(r,u)$. In the opposite direction,
a standard  estimate of the kernel 
$$
H(z) \ll \frac{|z|^2}{1+|z|}\,, \qquad z\in{\bf C}\,,
$$
yields Borel's estimate 
\begin{equation}
M(r,u) \ll 
r\int_0^r \frac{\mu(t)}{t^2}\,dt 
+ r^2 \int_r^\infty \frac{\mu(t)}{t^3}\,dt \,. 
\label{2.borel}
\end{equation}
In particular, 
\begin{equation}
M(r,u) = \left\{
\begin{array}{ll}
o(r), & r\to 0\, \\ \nonumber \\
o(r^2), & r\to\infty\,.
\end{array} 
\right.
\label{1.M}
\end{equation}

\begin{procl}{Theorem~1}
Let $u(z)$ be a canonical integral (\ref{1.1}) of genus one, then
\begin{equation}
M(r,u) \ll r\int_0^r \frac{\n(t)}{t^2}\,dt + r^2 \int_r^\infty
\frac{\n(t)}{t^3}\, dt 
+ r^2 \int_r^\infty \frac{\mathfrak m (t,u)}{t^2}\,dt \,,
\label{1.3}
\end{equation}
where
$$
\mathfrak m (r,u) = 
\frac{1}{2\pi} \int_0^{2\pi} u^-(re^{i\theta}|\sin\theta)|)
\frac{d\theta}{\sin^2\theta}
$$is the Tsuji proximity function of $u$.
\end{procl}

If the function $u$ is non-negative in $\bf C$, then the proximity
function $\mathfrak m(r,u)$ vanishes, 
and applying Jensen's inequality we arrive at

\medskip\par\noindent{\bf Corollary 1. }{\em Let $u$ be a canonical
integral of genus one which is non-negative in $\bf C$. Then}
$$
\mu (r) \ll r\int_0^r \frac{\n(t)}{t^2}\,dt + r^2 \int_r^\infty
\frac{\n(t)}{t^3}\, dt\,.
$$

\medskip
As we explained in the introduction, this result immediately yields the
classical Marcinkiewicz inequality (\ref{0.1}). 
In this case one can apriori assume that the function $f$ is
bounded, so that most of the technicalities needed for the proof of
Theorem~1 (see Lemmas~2 and 4 below) are redundant, and our proof of
Marcinkiewicz' inequality, being conceptually new, is comparable in
length to the classical one. 

There is a curious reformulation of Corollary~1. 
Let $\mathcal M$ be a measurable space endowed with a locally finite
non-negative measure $dm$, let $f$ be a complex-valued measurable
function on $\mathcal M$, and let 
$m_f(\lambda)=m(\{|f|\ge \lambda\})$ be the distribution function of $f$.
If
\begin{equation}
\int_{\mathcal M} \min(|f|,|f|^2)\, dm <\infty\,,
\label{?}
\end{equation}
then we define the logarithmic determinant of $f$ of genus one
$$
u_f(z) = \int_{\mathcal M} H(zf(t))\, dm(t)\,,
$$
which is subharmonic in $\bf C$, and moreover is represented by
a canonical integral of genus one. Applying Corollary~1, we obtain

\medskip\par\noindent{\bf Corollary~2. }{\em If $f$ satisfies
condition (\ref{?}), and the
logarithmic determinant $u_f$ is non-negative in $\bf C$, then}
$$
m_f(\lambda) \ll 
\left\{
\frac{1}{\lambda^2}\, \int_0^\lambda sm_{{\rm Im}\,f}(s) ds 
+ \frac{1}{\lambda}\, \int_\lambda^\infty m_{{\rm Im}\,f}(s) ds 
\right\}\,, \qquad 0<\lambda<\infty\,.
$$
{\em In particular, 
$$
||f||_{L^p(m)} \ll_p ||{\rm Im}\, f||_{L^p(m)}\,, \qquad 1<p<2\,,
$$
and}
$$
m_f(\lambda) \ll \frac{||{\rm Im}\, f||_{L^1(m)}}{\lambda}\,, \qquad
0<\lambda<\infty\,.
$$

\medskip
Corollary~1  can also be applied  to Jensen measures in $\bf C$. A
compactly supported finite measure $\sigma$ in $\bf C$ is called a
{\em Jensen measure} (with respect to the origin) if for an arbitrary
subharmonic function $h$ in $\bf C$
\begin{equation}
h(0) \le \int h d\sigma\,.
\label{10}
\end{equation}
A simple argument shows that (\ref{10}) then holds true for subharmonic
functions in a domain $G$ such that $0\in G$ and 
$\rm{supp}(\sigma)\subset G$. For a harmonic function, the equality sign
must occur in (\ref{10}). Therefore,
$$
\sigma ({\bf C}) = \int 1 d\sigma =1\,,
$$
that is, $\sigma $ is a probability measure, and 
\begin{equation}
\int \zeta^k d\sigma (\zeta) = 0\,, \qquad k=1,2,...\,.
\label{20}
\end{equation}

Define the potential
\begin{equation} 
V_{\sigma} (z) = \int \log|1-z\zeta| d\sigma (\zeta)\,.
\label{30}
\end{equation}
Then, due to (\ref{10}) and (\ref{20}), 
\begin{equation}
0\le V_\sigma(z) \le \log^+(c|z|)\,, \qquad z\in{\bf C}\,,
\label{40}
\end{equation}
for some $c>0$. The opposite is also true: if, for some $c>0$, a
subharmonic function $V$ satisfies (\ref{40}), then it is a potential of a
of a Jensen measure \cite[\S 4]{khabib}. 

Due to condition (\ref{20}), every potential $V_\sigma$ of a Jensen
measure can be represented by a canonical integral of genus one:
$$
V_\sigma(z) = \int_{\bf C} H(z/\zeta) d\mu(\zeta)\,, \qquad
d\mu(\zeta) = d\sigma(1/\zeta)\,.
$$
Thus, Theorem 1 is applicable to the potential $V_\sigma$, and we obtain

\medskip\par\noindent{\bf Corollary 3. }{\em Let $\sigma$ be a Jensen
measure in $\bf C$, $\sigma(\lambda) = \sigma(|z|\ge\lambda)$,
$\sigma_I(\lambda) = \sigma(|\rm{Im} z|\ge \lambda)$. Then}
$$
\sigma (\lambda) \ll \frac{1}{\lambda^2} 
\int_0^\lambda s\sigma_I(s) ds 
+ \frac{1}{\lambda} \int_{\lambda}^\infty \sigma_I(s) ds\,.
$$

\medskip The class of Jensen measures is invariant with respect to the
holomorphic mappings. More precisely, let $G$ be a domain which
contains the origin and $\rm{supp}(\sigma)$, and let $F$ be an analytic
function in $G$, $F(0)=0$. Then the push forward $F_* \sigma$ is defined
by
$$
\int \phi \,dF_* \sigma = \int \phi \circ F\, d\sigma\,,
$$
where $\phi $ is an arbitrary continuous function in $\bf C$, and
$\phi\circ F$ is a composition of $\phi$ and $F$.  By the
monotone convergence theorem this equation also holds for semicontinuous
functions. The measure $F_*\sigma $ automatically has a compact support
since $F$ is bounded on $\rm{supp}(\sigma)$. If $h$ is subharmonic in $\bf
C$, then $h \circ F$ is subharmonic in $G$, and 
$$
\int h\, dF_*\sigma = \int h\circ F\, d\sigma \ge h(F(0)) = h(0)\,. 
$$
Hence, $F_*\sigma$ is a Jensen measure. 

\medskip\par\noindent{\bf Corollary 4. }{\em Let $\sigma$ be a Jensen
measure in $\bf C$, and let $f=g+ih$ be an analytic function in a domain
$G$ which contains the origin and $\rm{supp}(\sigma)$, and  
$f(0)=0$. Let 
$$
m_{f,\sigma}(\lambda) = \sigma(|f|\ge \lambda)\,, 
\qquad m_{h,\sigma}(\lambda) = \sigma(|h|\ge \lambda)\,.
$$
Then}
\begin{equation}
m_{f,\sigma}(\lambda) \ll \frac{1}{\lambda^2} 
\int_0^\lambda sm_{h,\sigma}(s) ds 
+ \frac{1}{\lambda} \int_{\lambda}^\infty m_{h,\sigma}(s) ds\,.
\label{50}
\end{equation} 

\medskip Corollary~4 probably holds true under a weaker (and
more natural) assumption $g(0)=0$ rather than $f(0)=0$. In that case, 
using Theorem~2 (see below) one can get an estimate which is slightly
weaker than (\ref{50}). 

\medskip
In the next result,
we shall not assume that $u(z)$ is non-negative in $\bf C$ and instead
introduce
the quantity
$$
\delta (r) = \n(r) + [u^-(r) + u^-(-r)]
$$
which we keep under control. 

We assume that the integrals 
\begin{equation}
\int_0 \frac{\delta (t)}{t^2}\, dt\,, \qquad \mbox{and} \qquad
\int^\infty \frac{\delta (t)}{t^3} (1+\log t)\, dt
\label{1.a}
\end{equation}
are convergent and define 
\begin{equation}
\delta^* (r) = r \int_0^r \frac{\delta(t)}{t^2}\, dt + 
r^2 \int_r^\infty \frac{\delta (t)}{t^3} 
\left(1+\log\frac{t}{r} \right)\,dt
\label{1.*}
\end{equation}
The function $\delta^*(r)$ does not decrease, $r^{-2}\delta^*(r)$ does not
increase, 
and therefore
\begin{equation}
\delta^*(r) \le \delta^*(2r) \le 4\delta^* (r)\,, \qquad 0<r<\infty\,.
\label{1.d}
\end{equation} 

\begin{procl}{Theorem~2} Let $u(z)$ be an arbitrary subharmonic function
in $\bf C$ represented by a canonical integral of genus one. Then  
\begin{equation}
M(r,u) \ll r^2 
\left[ \int_r^\infty \frac{\sqrt{\delta^*(t)}}{t^2}\,dt \right]^2\,.
\label{1.4}
\end{equation}
\end{procl}

Observe, that the RHSs of (\ref{1.3}) and (\ref{1.4}) do not depend on
the bound for the integral (\ref{1.2}). Estimate (\ref{1.4}) is slightly
weaker than (\ref{1.3}); however, it suffices for 
deriving estimates of M.~Riesz and Kolmogorov, as well as of the weak
$(p,\infty)$-type estimate (see Corollary~6 below).

Fix an arbitrary  $\epsilon>0$. 
Then by the Cauchy inequality
\begin{eqnarray*}
\left[ \int_r^\infty \frac{\sqrt{\delta^*(t)}}{t^2}dt \right]^2 
&=& 
\left[ \int_r^\infty 
\frac{\sqrt{\left(1 +  
\log^{1+\epsilon}\frac{t}{r}\right)\delta^*(t)}}{t^{3/2}}\,
\frac{dt}{t^{1/2}\sqrt{1+\log^{1+\epsilon}\frac{t}{r}}}
\right]^2 \\ \\ 
&\ll_{\epsilon}& \int_r^\infty \frac{\delta^*(t)}{t^3}
\left(1+\log^{1+\epsilon}\frac{t}{r} \right)\,dt \\ \\ 
&\ll_{\epsilon}& 
\frac{1}{r} \int_0^r \frac{\delta(s)}{s^2}\,ds +  
\int_r^\infty \frac{\delta(s)}{s^3} \left(1+\log^{3+\epsilon}\frac{s}{r} 
\right)\, ds\,.
\end{eqnarray*}
Thus we get

\begin{procl}{Corollary 5} For each $\epsilon>0$,
\begin{equation}
M(r,u) \ll_{\epsilon} 
r \int_0^r \frac{\delta(t)}{t^2}\,dt + r^2 
\int_r^\infty \frac{\delta(t)}{t^3} \left(1+\log^{3+\epsilon}\frac{t}{r} 
\right)\, dt\,.
\label{1.4a}
\end{equation}
\end{procl}

We do not know whether the term $\log^{3+\epsilon}$ is really needed on
the RHS of (\ref{1.4a}). 
Apparently, our method does not allow us to omit it. 
Rewriting (\ref{1.4a}) in the form
$$
M(r,u) \ll_{\epsilon}
\int_0^1 \frac{\delta (rs)}{s^2}\, ds 
+ \int_1^\infty \frac{\delta (rs)}{s^3} 
\left( 1+ \log^{3+\epsilon}s \right)\, ds\,, 
$$
we immediately obtain 

\begin{procl}{Corollary 6} The following inequalities hold for canonical
integrals of genus one:

\par\noindent M. Riesz-type estimate:
\begin{equation}
\int_0^\infty \frac{\mu(r)}{r^{p+1}}\, dr
\ll_p \int_0^\infty \frac{M(r,u)}{r^{p+1}}\, dr
\ll_p\int_0^\infty \frac{\delta(r)}{r^{p+1}}\, dr\,,
\qquad 1<p<2\,, 
\label{1.5}
\end{equation}
weak $(p,\infty)$-type estimate: 
\begin{equation}
\sup_{r\in (0,\infty)} \frac{\mu(r)}{r^p} \ll_p
\sup_{r\in (0,\infty)} \frac{M(r,u)}{r^p} \ll_p
\sup_{r\in (0,\infty)} \frac{\delta (r)}{r^p}\,,
\qquad 1<p<2\,,
\label{1.6}
\end{equation}
and
Kolmogorov-type estimate: 
\begin{equation}
\sup_{r\in (0,\infty)} \frac{\mu(r)}{r} \ll
\sup_{r\in (0,\infty)} \frac{M(r,u)}{r} \ll
\int_0^\infty \frac{\delta(r)}{r^2}\, dr\,.
\label{1.7}
\end{equation}
\end{procl}

Estimates weaker than (\ref{1.5}) and (\ref{1.7}) were obtained in
\cite{MS1} and \cite{MS3} under additional restrictions which now appear
to be redundant. Estimate (\ref{1.6}) is apparently new.  

If we assume that $d\mu$ is supported by $\bf R$, that is, $u(z)$ is
harmonic in $\bf C_{\pm}$, then our technique gives a better result:

\smallskip
{\em Let $u(z)$ be a canonical integral of genus one
of a measure $d\mu $ supported by $\bf R$. Then, for} $0<r<\infty$,
\begin{equation}
M(r,u) \ll
r \int_0^r \frac{u^-(t)+u^-(-t)}{t^2}\, dt
+ r^2 \int_r^\infty \frac{u^-(t)+u^-(-t)}{t^3}
\, dt\,.
\label{1.8}
\end{equation}

Note, that one cannot replace $u^-(x)$ by 
$u^+(x)$ on the RHS of our estimates. For example, the function
$$
u_p(z) = r^p \cos p \left(\frac{\pi}{2} - |\theta| \right),
\qquad 1<p<2\,,
$$
is subharmonic in $\bf C$, harmonic in the upper and 
lower half-planes $\bf C_\pm$, represented by the
canonical integral of genus one of the measure 
$d\mu (x) = c_p |x|^{p-1}dx$ ($c_p>0$), and non-positive on $\bf R$.
 
There is a corollary to Theorem~2 which is parallel to Corollary~2. 
Let $\mathcal M$ be a measurable space endowed with a locally
finite non-negative measure $dm$, and let 
$f:\, {\mathcal M}\to {\bf R}^{n+1}$, $n\ge 1$, be a measurable function
such that
\begin{equation}
\int_{\mathcal M} \min(||f||, ||f||^2)\, dt <\infty\,,
\label{1.9a}
\end{equation}
where $||\,.\,||$ stands for the $n+1$-dimensional Euclidean norm.
We start to enumerate the coordinates in ${\bf R}^{n+1}$
with $j=0$, and denote by $e_0$ the vector in ${\bf R}^{n+1}$ with the
zeroth coordinate equal one, and other coordinates vanishing.
Let $f_j(t)$ be the $j$-th coordinate function of $f(t)$,
and ${\hat f}(t) =\left\{\sum_{j=1}^n f_j^2(t) \right\}^{1/2}$.
We define the
logarithmic determinant
\begin{eqnarray*}
v_f(x) &=& \int_{\mathcal M}
\big[ \log||e_0-xf(t)||+xf_0(t)  \big]\, dm(t) \\ \\
&=& \int_{\mathcal M} 
\left[\log\sqrt{1-2xf_0(t)+x^2||f||^2} + xf_0(t)\right]\, dm(t)\,,
\qquad x\in {\bf R}\,,
\end{eqnarray*}
where the integral converges due to assumption (\ref{1.9a}). Then, if the
function $v_f(x)$ is non-negative on $\bf R$, we may estimate its
distribution function $m_f(\lambda)=m(\{||f||\ge \lambda \})$ by the
distribution function $m_{\hat f}=m(\{\hat f \ge \lambda
\})$ of $\hat f$. 

For this, observe that 
$$
v_f(x) = \int_{\mathcal M} H(xf_C(t))\, dm(t)\,,
\qquad x\in {\bf R}\,,
$$
where $f_C$ is a ``complex-valued surrogate''  of $f$:
$f_C=f_0+i\hat{f}$. That is, $v_f$ has a subharmonic continuation 
from $\bf R$ to $\bf C$ by a canonical integral of genus one
$$
u_{f_C}(z) = \int_{\mathcal M} H(zf_C (t))\, dm(t)\,, 
\qquad z\in {\bf C}\,.
$$
Next, observe that 
$m_f(\lambda) = m(\{f_0^2+{\hat f}^2\ge \lambda^2\})=m_{f_C}(\lambda)$,
and $m_{\hat f}(\lambda) = m_{{\rm Im}\, f_C}(\lambda)$ 
for $0<\lambda<\infty$. Hence, Theorem~2 is applicable in this situation.
For simplicity, we restrict
ourselves to the case when $v_f$ is non-negative on the real axis.

\medskip\par\noindent{\bf Corollary~7. }{\em Let $f$ satisfy 
condition (\ref{1.9a}), and let the logarithmic determinant $v_f$ be
non-negative on
the real axis. Then, for $0<\lambda<\infty$ and $\epsilon>0$,}
$$
m_f(\lambda) \ll_{\epsilon}
\frac{1}{\lambda^2} \int_0^\lambda s
\left(1+\log^{3+\epsilon}\frac{\lambda}{s}\right)m_{\hat f}(s)\,ds
+\frac{1}{\lambda} \int_\lambda^\infty m_{\hat f}(s)\,ds\,.
$$  
{\em In particular,} 
$$
||f||_{L^p(m)} \ll_p ||\hat f||_{L^p(m)}\,, \qquad 1<p<2\,,
$$
{\em and}
$$
m_f(\lambda) \ll \frac{||\hat f||_{L^1(m)}}{\lambda}\,.
$$ 

This corollary may be of some interest in view of the results of 
Aleksandrov and Kargaev \cite{AK}. 

\medskip

Our third result pertains to a more general class of subharmonic
functions represented by a {\em generalized} canonical integral of
genus one. It gives a Kolmogorov-type estimate which 
can be applied to a wider class of functions than 
(\ref{1.7}): 
 
\begin{procl}{Theorem 3}
Let $d\mu$ be a non-negative locally finite measure on $\bf C$
such that
$$
\int_{\{|\zeta| \ge 1\}} \frac{d\mu(\zeta)}{|\zeta|^{2}} < \infty\,,
$$
and let there exist a finite principal value integral
$$
\lim_{\varepsilon\to 0}\,
\int_{\{\varepsilon \le |\zeta|\le 1\}}
\frac{d\mu(\zeta)}{\zeta}\,.
$$
Let 
$$
u(z)=\lim\limits_{\varepsilon\to 0} \,\int_{|\zeta|>\varepsilon}
H(z/\zeta)\, d\mu(\zeta)\,,
$$
then
\begin{equation}
\sup_{0<r<\infty} \frac{M(r,u)}{r} \ll
\int_0^\infty \frac{\delta (t)}{t^2}\, dt +
\limsup\limits_{r\to 0}\, \frac{\mu(r)}{r}\,.
\label{kolmogorov}
\end{equation}
\end{procl}

It is easy to see that if the integral (\ref{1.2}) converges at the
origin, then the upper limit on the RHS of (\ref{kolmogorov}) vanishes,
and in this case (\ref{kolmogorov}) coincides with (\ref{1.7}).

In fact, our proof yields a stronger result 
\begin{equation}
\int_{-\infty}^\infty \frac{u^+(t)}{t^2}\, dt + 
\limsup_{r\to\infty} \frac{M(r,u)}{r} 
\ll \int_0^\infty \frac{\delta (t)}{t^2}\, dt +
\limsup\limits_{r\to 0}\, \frac{\mu(r)}{r}\,,
\label{cart}
\end{equation} 
which gives control over the positive harmonic majorants of $u$ in the
upper and lower half-planes. Applying a known technique of functions of
Cartwright class \cite{Levin}, \cite{Koosis1}, one can extract from
(\ref{cart}) information about the asymptotic regularity of $u$ and
$\mu$ at infinity and near the origin.

Notice, that one can reformulate Theorem~3 in the spirit of Corollaries~2
and 7. We leave this to the reader.

\section{Auxiliary Lemmas}
\setcounter{equation}{0}

We shall need several known facts about harmonic and subharmonic
functions.

\begin{procl}{Lemma 1} 
Let $v$ be a  subharmonic function in the angle $S=\{z:\, 0<\arg z <
\alpha\}$, 
$0<\alpha<2\pi$, let
\begin{equation}
\limsup_{z\to\zeta, \, z\in S} v^+(z) \le \Phi(|\zeta|)\,, \qquad
\zeta\in \partial S\,; 
\label{2.a}
\end{equation}
and let 
\begin{equation}
\int_0^\alpha v^+(re^{i\theta}) \sin
\left(\frac{\pi}{\alpha}\theta\right) \, d\theta
= o(r^{\pi/\alpha})\,, \qquad r\to\infty\,.
\label{2.b}
\end{equation}
Then, for $z=re^{i\theta}\in S$,
\begin{equation}
v(re^{i\theta}) \sin\left(\frac{\pi}{\alpha}\theta\right) \ll_{\alpha}  
r^{-\pi/\alpha} \int_0^r \Phi (t) t^{\pi /\alpha -1}\, dt 
+ r^{\pi/\alpha} \int_r^\infty \frac{\Phi (t)}{t^{\pi /\alpha +1}}\, dt\,. 
\label{2.4}
\end{equation}
If the majorant $\Phi (t)$ does not decrease, then the factor
$\sin(\pi \theta/\alpha)$ on the LHS of (\ref{2.4}) can be omitted. 
\end{procl}

\par\noindent{\em Proof:} The general case is easily reduced to the 
special case when $S=\bf C_+$, so that, without loss of generality, 
we assume that $\alpha=\pi$. First, we show that $v(z)$ is majorized 
by the Poisson integral of $\Phi (|t|)$, and 
then we estimate this integral.

Denote by $h_R(z)$ a harmonic function in the semi-disk
$\{\mbox{Im}z>0, |z|<R\}$ with boundary values $h_R(t)=\Phi(|t|)$,
$-R<t<R$, and $h_R(Re^{i\theta})=v^+(Re^{i\theta})$, $0<\theta<\pi$. 
Applying the Poisson-Nevanlinna representation in this semi-disk
(see \cite[Chapter~1, Theorem~2.3]{GO}, \cite[Section~24.3]{Levin}), we
obtain for $z=re^{i\theta}$, $r<R$,
\begin{equation}
v(z)\le h_R(z) 
= \int_{-R}^R \Phi (|t|) K_1(z,t)\,dt 
+ \int_0^\pi v^+(Re^{i\phi}) 
K_2(z, Re^{i\phi})\, d\phi\,, 
\label{2.c}
\end{equation} 
where
\begin{equation}
K_1(z,t) = \frac{r\sin\theta}{\pi}
\left\{ \frac{1}{|z-t|^2} - \frac{R^2}{|R^2-zt|^2}\right\}\,,
\label{*1}
\end{equation}
\begin{equation}
K_2(z,Re^{i\phi})=\frac{1}{2\pi} \frac{4Rr(R^2-r^2)\sin\phi\sin\theta}
{(R^2+r^2-2Rr\cos(\phi-\theta))(R^2+r^2-2Rr\cos(\phi+\theta))}\,,
\label{*2}
\end{equation}

By condition (\ref{2.b}), the second integral on the RHS of 
(\ref{2.c}) tends to $0$ as $R\to\infty$. Therefore, letting 
$R\to\infty$ in (\ref{2.c}), we obtain
\begin{equation}
v(z) \le \frac{r\sin\theta}{\pi}
\int_{-\infty}^{\infty} \frac{\Phi (|t|)}{|z-t|^2}\, dt\,.
\label{2.d}
\end{equation}
Making use of straightforward estimates of the Poisson 
kernel, we get
$$
v(z) \le \frac{2}{\pi r\sin\theta} \int_0^{2r} \Phi (t)\, dt
+ \frac{4r}{\pi} \int_{2r}^\infty \frac{\Phi(t)}{t^2}\, dt\,,
$$
and estimate (\ref{2.4}) follows.

If the majorant $\Phi (t)$ does not decrease, then we modify the
previous argument:
$$
v(z) \le \frac{4}{r} \int_0^{r/2} \Phi (t)\, dt + \Phi (2r)
+ 4r \int_{2r}^\infty \frac{\Phi (t)}{t^2}\, dt 
\ll \frac{1}{r} \int_0^r \Phi (t)\, dt + 
+ r \int_r^\infty \frac{\Phi (t)}{t^2}\, dt\,,
$$
completing the proof. $\Box$ 

\medskip

The next lemma asserts that under certain conditions the Carleman
integral formula \cite[Lecture~24]{Levin}, \cite[Chapter~1]{GO}
holds without remainder.

\begin{procl}{Lemma 2}
Let $v(z)$ be a subharmonic function on $D_R=\{z\in {\bf \bar C_+}:\,
|z|\le
R\}$
which satisfies conditions
\begin{equation}
\int_0^\pi v^+(re^{i\theta}) \sin\theta\, d\theta = o(r)\,,
\qquad r\to 0\,,
\label{i}
\end{equation}
and
\begin{equation}
\int_0 \frac{\delta (t)}{t^2}\, dt <\infty\,.
\label{ii}
\end{equation}
Then 
\begin{eqnarray}
\frac{1}{2\pi} \int_{-R}^R v(t) 
\left( \frac{1}{t^2} - \frac{1}{R^2} \right)dt 
&+&\frac{1}{\pi R} \int_0^\pi v(Re^{i\phi}) \sin\phi\, d\phi 
\nonumber \\ \nonumber \\
&=&
\int_{D_R}
\left( \frac{1}{|\zeta|^2}-\frac{1}{R^2}\right)
{\rm Im}\zeta
\, d\mu(\zeta)\,,
\label{2.6}
\end{eqnarray}
where the first integral on the LHS is absolutely convergent. 
\end{procl}

\par\noindent{\em Proof:} 
We start with the Nevanlinna representation
\begin{eqnarray*}
\int_{-R}^R v^+(t) K_1(z,t)\,dt 
&+& \int_0^\pi v^+(Re^{i\phi}) K_2(z,Re^{i\phi})\, d\phi \\ \\
&=& v(z) +  
\int_{-R}^R v^-(t) K_1(z,t)\,dt 
+ \int_0^\pi v^-(Re^{i\phi}) K_2(z,Re^{i\phi})\, d\phi \\ \\
&\ & \qquad +  \int_{D_R} K_3(z,\zeta)\, d\mu(\zeta)\,,
\end{eqnarray*}
where the kernels $K_1$ and $K_2$ were defined by (\ref{*1}) and
(\ref{*2}), and
\begin{equation}
K_3(z,\zeta) = \log\left|
\frac{z-\bar\zeta}{z-\zeta}\,\cdot\,
\frac{R^2-z\bar\zeta}{R^2-z\zeta}\right|\,.
\label{*3}
\end{equation}

We multiply both the left and right hand sides of the Nevanlinna
representation by 
$r^{-1}\sin\theta$, integrate it with respect to $\theta$ from $0$ to
$\pi$
and change the integration order in all terms. We shall use the 
formulas:
\begin{eqnarray}
\frac{1}{r} 
\int_0^\pi K_1(re^{i\theta},t) \sin\theta\,d\theta = \frac{1}{2}
\left[\min\left(\frac{1}{t^2}-\frac{1}{r^2} \right)-\frac{1}{R^2}
\right]\,,
\label{*4}
\end{eqnarray}
\begin{eqnarray}
\frac{1}{r} \int_0^\pi K_2(re^{i\theta}, Re^{i\phi}) \sin\theta\, d\theta  
=\frac{1}{R} \sin\phi\,,
\label{*5}
\end{eqnarray}
and
\begin{eqnarray}
\frac{1}{r} \int_0^\pi K_3(re^{i\theta},\zeta) \sin\theta\, d\theta
= \pi \mbox{Im}\zeta\,\left[
\min\left(\frac{1}{|\zeta|^2},\frac{1}{r^2}\right)-\frac{1}{R^2}
\right] \,.
\label{*6}
\end{eqnarray}

Observe that the RHS of relations (\ref{*4})-(\ref{*6}) are non-decreasing
functions of $r^{-1}$. Therefore, making the limit transition $r\to 0$,
and using the monotone convergence theorem and condition (\ref{ii}) of the
lemma, we get
\begin{eqnarray*}
\frac{1}{2} \int_{-R}^R v^+(t) 
\left( \frac{1}{t^2} - \frac{1}{R^2} \right)dt 
&+& \frac{1}{R} \int_0^\pi v^+(Re^{i\phi}) \sin\phi d\phi 
\nonumber \\ \nonumber \\
&=&\frac{1}{2} \int_{-R}^R v^-(t) 
\left( \frac{1}{t^2} - \frac{1}{R^2} \right)dt 
+ \frac{1}{R} \int_0^\pi v^-(Re^{i\phi}) \sin\phi\ d\phi \\ \\  
&\,& \qquad + \pi \int_{D_R}
\left( \frac{1}{|\zeta|^2}-\frac{1}{R^2}\right)
\mbox{Im}\zeta
\, d\mu(\zeta)\,.
\end{eqnarray*}
The first and third integrals  on the RHS are finite due to condition
(\ref{ii}). This completes the proof. $\Box$

\begin{procl}{Remark}
Condition (\ref{i}) holds true for canonical integrals of genus one
defined in (\ref{1.1}).
\end{procl}
Indeed, if $u(z)$ is such an integral, then due to (\ref{1.3})
$$
\int_0^{2\pi} u^+(re^{i\theta})\, d\theta = o(r), \qquad r\to 0\,.
$$
Since $u(0)=0$, this yields
\begin{eqnarray*}
\int_0^{2\pi} |u(re^{i\theta})|\, d\theta 
&=& 2 \int_0^{2\pi} u^+(re^{i\theta})\, d\theta
- \int_0^{2\pi} u(re^{i\theta})\, d\theta \\ \\ 
&=& 
2 \int_0^{2\pi} u^+(re^{i\theta})\, d\theta
= o(r)\,, \qquad r\to 0\,.
\end{eqnarray*}

The third lemma was proved in \cite{MS3} (cf. \cite[Lecture~26]{Levin}). 
Its proof uses the Nevanlinna representation for the semi-disk.
 
\begin{procl}{Lemma 3}
Let $v(z)$  be a function which is harmonic in $\bf C_+$,
subharmonic in $\overline{\bf C}_+$, and satisfies conditions
(\ref{i}) and (\ref{ii}) of Lemma~2. Then,
for $z=re^{i\theta}\in {\bf C_+}$,
\begin{equation}
v(re^{i\theta}) \sin\theta 
\ll \frac{1}{\pi} \int_0^\pi v^-(2re^{i\varphi}) \sin\varphi \, d\varphi
+ \frac{r}{2\pi} \int_{-2r}^{2r} \frac{v^-(t)}{t^2}\, dt\,.
\label{2.7}
\end{equation}
\end{procl}

The next lemma is a version of the Levin integral formula
(\cite{Levin1}, \cite[Chapter~1]{GO}) without
a remainder. 

\begin{procl}{Lemma~4} Let $v$ be a subharmonic function in $\bf C$ such
that $v(z)$ and $v(\bar z)$ satisfy conditions (\ref{i}) and (\ref{ii}) of
Lemma~2. Then
\begin{equation}
\frac{1}{2\pi}
\int_0^{2\pi} v(Re^{i\theta}|\sin\theta|) \frac{d\theta}{R\sin^2\theta}
= \int_0^R \frac{\n(t)}{t^2}\, dt\,,
\label{*7}
\end{equation}
where $\n(t)$ is the Levin-Tsuji counting function, and the integral
on the LHS is absolutely convergent.
\end{procl}

\par\noindent{\em Proof:} 
It suffices to prove that
\begin{equation}
\frac{1}{2\pi} \int_0^\pi v(Re^{i\theta}\sin\theta)
\frac{d\theta}{R\sin^2\theta}
=\int_{\left|{\rm Im}\frac{1}{\zeta}\right|>\frac{1}{R}}
\left[\left|{\rm Im} \frac{1}{\zeta}\right|-\frac{1}{R}
\right]\,d\mu(\zeta)\,.
\label{??}
\end{equation}
Then (\ref{*7}) follows by adding to (\ref{??}) a similar formula for the
integral from $\pi$ to $2\pi$.

First, we prove that the integral on the LHS of
relation (\ref{*7}) is absolutely convergent. 
Making use of notations introduced in (\ref{*1}), (\ref*2) and (\ref{*3}),
observe that the Nevanlinna formula implies that
$$
|v(z)| \le 
\int_{-R}^R |v(t)|K_1(z,t)\, dt + \int_0^\pi |v(Re^{i\phi})|
K_2(z,Re^{i\phi}) \,d\phi + \int_{D_R} K_3(z,\zeta)\,
d\mu(\zeta)\,.
$$
We set $z=Re^{i\theta}\sin\theta$, multiply the formula by
$(R\sin^2\theta)^{-1}$, integrate it with respect to $\phi$ from $0$ to
$\pi$, and change the integration order in all terms. We shall use the
following relations:
$$
\int_0^\pi K_1(Re^{i\theta}\sin\theta, t) \frac{d\theta}{R\sin^2\theta}
= \frac{1}{t^2}-\frac{1}{R^2}\,,
$$
$$
\int_0^\pi K_2(Re^{i\theta}\sin\theta, Re^{i\phi})
\frac{d\theta}{R\sin^2\theta}
= \frac{2}{R}\sin\phi\,,
$$
and 
$$
\int_0^\pi K_3(Re^{i\theta}\sin\theta, \zeta)
\frac{d\theta}{R\sin^2\theta}
= 2\pi \left[ \min
\left(\left|\mbox{Im}\frac{1}{\zeta}\right|,\frac{1}{R}\right)
- \frac{\mbox{Im}\zeta}{R^2}
\right]\,.
$$ 
Using these relations, we verify that 
\begin{eqnarray*}
\int_0^\pi |v(Re^{i\theta}\sin\theta)|\, \frac{d\theta}{R\sin^2\theta}
&\le& \int_{-R}^R |v(t)|
\left(\frac{1}{t^2}-\frac{1}{R^2}\right)\, dt \\ \\
&+& \frac{2}{R} \int_0^\pi |v(Re^{i\phi})|\sin\phi\,d\phi
+2\pi \int_{D_R}\left|\mbox{Im}\frac{1}{\zeta}\right|\, d\mu(\zeta)\,.
\end{eqnarray*}
The first integral on the RHS is finite due to Lemma~2, and the third is
finite due to condition (\ref{ii}). That is, the integral on the LHS of
(\ref{*7}) is absolutely convergent. 

Now, we write the Nevanlinna formula in the form
$$
v(z)=\int_{-R}^R v(t)K_1(z,t)\, dt + \int_0^\pi
v(Re^{i\phi})K_2(z,Re^{i\phi})\,d\phi 
- \int_{D_R} K_3(z,\zeta)\, d\mu(\zeta)\,.
$$
Again, we set here $z=Re^{i\theta}\sin\theta$, multiply by  
$(R\sin^2\theta)^{-1}$, integrate with respect to $\theta$ 
from $0$ to $\pi$ and change the integration order in all terms. 
We can do this since we already know that the integrals with 
$|v|$ instead of $v$ are finite. As a result, we obtain the equation
\begin{eqnarray}
\int_0^\pi v(Re^{i\theta}\sin\theta)\frac{d\theta}{R\sin^2\theta} 
&=& \int_{-R}^R v(t) \left(\frac{1}{t^2}-\frac{1}{R^2} \right)\, dt
+ \frac{2}{R} \int_0^\pi v(Re^{i\phi})\sin\phi\, d\phi 
\nonumber \\ \nonumber \\
&\,& \ - 2\pi \int_{D_R} \left[ 
\min\left(\left|\mbox{Im}\frac{1}{\zeta}\right|,\frac{1}{R}\right) - 
\frac{\mbox{Im}\zeta}{R^2}\right]\, d\mu(\zeta)\,. 
\nonumber \\
\label{*8}
\end{eqnarray}
Taking into account (\ref{2.6}), we get
\begin{eqnarray*}
\int_0^\pi v(Re^{i\phi}\sin\phi)\frac{d\phi}{R\sin^2\phi} 
&=& 2\pi \int_{D_R} \left[ \frac{1}{|\zeta|^2} - \frac{1}{R^2} \right]
\mbox{Im}\zeta \, d\mu(\zeta) \\ \\ 
&\,&\qquad - 2\pi \int_{D_R} \left[ 
\min\left(\left|\mbox{Im}\frac{1}{\zeta}\right|,\frac{1}{R}\right) - 
\frac{\mbox{Im}\zeta}{R^2}\right]\, d\mu(\zeta)
\\ \\
&=&
2\pi \int_{|{\rm Im}\frac{1}{\zeta}|\ge\frac{1}{R}} \left[ 
\left|\mbox{Im}\frac{1}{\zeta}\right|
- \frac{1}{R} \right]\, d\mu(\zeta)\,. 
\end{eqnarray*}
Then (\ref{??}) follows and the proof is complete.
$\Box$

\medskip In other words, in the assumptions of Lemma~4, 
the first fundamental theorem for Tsuji characteristics holds
without a remainder term:
\begin{equation}
\T(r,u) = \mathfrak m(r,u) + \int_0^r\frac{\n(t)}{t^2}\,dt\,,
\qquad 0<r<\infty\,,
\label{tsuji}
\end{equation}
where 
$$
\T(r,u) = \frac{1}{2\pi} \int_0^{2\pi} u^+(re^{i\theta}|\sin\theta)|)
\frac{d\theta}{\sin^2\theta}\,,
$$
and 
$$
\mathfrak m (r,u) = 
\frac{1}{2\pi} \int_0^{2\pi} u^-(re^{i\theta}|\sin\theta)|)
\frac{d\theta}{\sin^2\theta}\,.
$$

\medskip
The last lemma was proved in a slightly different setting in \cite{LO}
(see also \cite[Lemma~5.2, Chapter~6]{GO}):

\begin{procl}{Lemma~5}
Let $u(z)$ be a subharmonic function in $\bf C$, and let 
$$
T(r,u) = \frac{1}{2\pi} \int_0^{2\pi} u^+(re^{i\theta})\,d\theta 
$$
be its Nevanlinna characteristic function.
Then, for $0<R<\infty$,
\begin{equation}
\int_R^\infty \frac{T(r,u)}{r^3}\, dr \le \int_R^\infty
\frac{\T(r,u)}{r^2}\, dr\,.
\label{tsuji1}
\end{equation}
\end{procl}

\section{Proof of Theorem 1}
\setcounter{equation}{0}

Using monotonicity of $T(r,u)$, Lemma~5, and then Lemma~4, we obtain
\begin{eqnarray*}
\frac{T(R,u)}{R^2} &\le& 2 \int_R^\infty \frac{T(r,u)}{r^3}\, dr \\ \\
&\stackrel{(\ref{tsuji1})}\le& 2 
\int_R^\infty \frac{\T(r,u)}{r^2}\, dr \\ \\ 
&\stackrel{(\ref{tsuji})}=& 2 \int_R^\infty \frac{dr}{r^2} 
\left( \int_0^r
\frac{\n(t)}{t^2}\, dt + \mathfrak m(r,u) \right)  \\ \\
&=& \frac{2}{R} \int_0^R \frac{\n(t)}{t^2}\, dt + 2\int_R^\infty 
\frac{\n(t)}{t^3}\, dt 
+ 2\int_R^\infty \frac{\mathfrak m(t,u)}{t^2}\,dt\,.
\end{eqnarray*}
Then the inequality $M(r,u)\le 3T(2r,u)$ completes the
proof. $\Box$

\section{Proof of Theorem 2}
\setcounter{equation}{0}

We split the proof into several parts. Without loss of generality, we 
assume convergence of the integrals
$$
\int_0 \frac{\delta(t)}{t^2}\, dt
\qquad \hbox{and} \qquad 
\int^\infty \frac{\delta(t)}{t^3} \log t\, dt\,.
$$

We define a measure $\mu_1$, ${\rm supp}(\mu_1)\subset {\bf{\bar C}_-}$, 
by reflecting at the real axis the part of the measure $\mu$ which lies in
the upper half-plane. Formally, 
$$
\mu_1(E) = \mu(E\cap {\bf\bar C}_-) + \mu(E^-\cap \bf{\bar C}_-)\,,
$$
where $E\subset \bf C$ is a borelian set, and $E^-=\{z:\, \bar z\in E\}$. 
Then the measure $\mu_1$ also satisfies condition (\ref{1.1}) and we
denote by $u_1(z)$ its canonical integral of genus one. 
Observe that $u_1(t)=u(t)$, so that $\delta(t,u_1) = \delta(t,u)$,
$t\in {\bf R}$.

\subsection{Estimate of $u_1^-(iy)$, $y>0$.}
We have
\begin{eqnarray*}
H(iy/\zeta) &=&  
\log \left| 1+ y\, \hbox{Im}\frac{1}{\zeta}
- iy\, \hbox{Re}\frac{1}{\zeta} \right| - 
y\, \hbox{Im}\frac{1}{\zeta} \\  \\
&\ge& \log\left| 1 + y\, \hbox{Im}\frac{1}{\zeta}\right| 
- y\, \hbox{Im}\frac{1}{\zeta}\,.
\end{eqnarray*}
Since the RHS is non-positive for $y>0$ and $\zeta\in {\bf {\bar C}_-}$,
\begin{eqnarray}
u_1^-(iy) &\le& -\int_{\bf{\bar C}_-} 
\left[\log\left| 1 + y\, \hbox{Im}\frac{1}{\zeta}\right| 
- y\, \hbox{Im}\frac{1}{\zeta}\right]\,d\mu(\zeta) \nonumber  \\ \nonumber \\
&=&-\int_0^\infty \left[ \log\left( 
1+\frac{y}{t} \right) - \frac{y}{t}\right] d\n(t) 
\nonumber \\ \nonumber \\
&=& y^2 \int_0^\infty \frac{\n(t)}{t^2(t+y)}\, dt 
\nonumber \\ \nonumber \\
&\le& y\int_0^y \frac{\n(t)}{t^2}\, dt 
+ y^2 \int_y^\infty \frac{\n(t)}{t^3}\, dt\,.
\label{3.1}
\end{eqnarray}

\subsection{Estimates of $u_1^+(re^{i\theta })$, $0<\theta<\pi$.}

Using harmonicity of the function $u_1$ in the upper half-plane, we
transform the lower bound for $u_1$ into the
upper bound. We shall show that
\begin{equation}
u_1^+(re^{i\theta}) \sin\theta \ll \delta^* (r)\, 
\qquad 0<r<\infty\,, \quad 0<\theta<\pi\,,
\label{3.4}
\end{equation}
where $\delta^*(r)$ is defined by (\ref{1.*}).

Consider the function $-u_1(z)$ and apply Lemma~1 to the angles 
$\{0<\arg z < \pi/2\}$ and $\{\pi/2<\arg z <\pi\}$ with
$$
\Phi (r) = [u_1^-(r) + u_1^-(-r)] + r\int_0^r \frac{\n(t)}{t^2}\, dt
+ r^2 \int_r^\infty \frac{\n(t)}{t^3}\, dt\,.
$$
Condition (\ref{2.a}) holds due to estimate (\ref{3.1}), and condition
(\ref{2.b}) holds due to estimate (\ref{1.M}) combined with Jensen's
inequality:
$$
\int_0^\pi u_1^-(re^{i\theta})\,d\theta 
\le \int_0^{2\pi} u_1^+(re^{i\theta})\, d\theta \le
M(r,u_1) = o(r^2)\,, \qquad r\to\infty\,.
$$
Therefore,
\begin{eqnarray}
-u_1(re^{i\theta}) |\sin 2\theta| &\ll&
\frac{1}{r^2}\int_0^r \Phi (t) t\, dt 
+ r^2 \int_r^\infty \frac{\Phi (t)}{t^3}\, dt
\nonumber \\ \nonumber \\ 
&\ll&
\frac{1}{r^2}\int_0^r [u_1^-(t)+u_1^-(-t)] t\, dt 
+ r^2 \int_r^\infty \frac{u_1^-(t)+u_1^-(-t)}{t^3}\, dt
\nonumber \\ \nonumber \\ 
&\ & \qquad +  
r \int_0^r \frac{\n(s)}{s^2}\, ds +r^2 \int_r^\infty \frac{\n(s)}{s^3}
\left(1+\log\frac{s}{r} \right)\, ds
\nonumber \\ \nonumber \\ 
&\ll& \delta^*(r)\,.
\label{4.3.1}
\end{eqnarray}

Observe that the factor $|\sin 2\theta|$ on the LHS of (\ref{4.3.1}) can
be replaced by $\sin\theta$. This follows from inspection of the proof of
Lemma~1 (since on the imaginary axis the function $-u(iy)$ has an
increasing majorant). Alternatively, one may again apply Lemma~1
to a small angle around the imaginary axis, say in
$\{|\theta-\pi/2|<\pi/8\}$. That is, we have 
\begin{equation}
-u_1(re^{i\theta})\sin\theta
\ll \delta^*(r)\,.
\label{4.3.2}
\end{equation}

Using Lemma~3 we obtain
$$
u_1^+(re^{i\theta}) \sin\theta \ll
\int_0^{\pi} u_1^-(2re^{i\phi}) \sin\phi \,d\phi + 
r\int_0^{2r} \frac{u_1^-(t)+u_1^-(-t)}{t^2}\, dt 
\ll \delta^*(r)\,,
$$
proving estimate (\ref{3.4}). 

\subsection{Estimate of $u^+(re^{i\theta})$, $\theta\ne 0,\pi$.}

Here we prove that, for an arbitrary $\eta >0$,
\begin{eqnarray}
u^+(re^{i\theta}) \ll
\frac{\delta^*(r)}{\eta |\sin\theta|} +
\frac{\eta r^2}{\sin^2\theta} \int_r^\infty
\frac{M(t,u)}{t^3}\, dt\,.  
\label{5.a}
\end{eqnarray}

For this, we shall need several upper bounds for
the difference
$$
D=D(z,\zeta)=H(z/\zeta)-H(z/{\bar \zeta})
=\log\left|\frac{1-z/\zeta}{1-z/\bar\zeta} \right| + 
{\rm Re} \left[z\left( \frac{1}{\zeta} 
- \frac{1}{\bar \zeta}\right) \right]\,,
$$
when $z,\zeta \in {\bf \bar C}_+$.

First,
\begin{equation}
D = \log\left|\frac{z-\zeta}{z-\bar\zeta} \right| 
+ 2{\rm Im}z\, \left|{\rm Im}\frac{1}{\zeta} \right|
\le 2|z|\, \left|{\rm Im}\frac{1}{\zeta} \right|
\,.
\label{**1}
\end{equation}
We shall use this estimate when 
$|z| \left|{\rm Im}\frac{1}{\zeta}  \right| \ge 1$.

Next,  let $t=|z|/|\zeta|$, $\theta=\arg(z)$, $\phi=\arg(\zeta)$.
Then
\begin{eqnarray}
D &=& \frac{1}{2} \log
\left[1-\frac{4t\sin\theta\sin\phi}{|1-te^{i(\theta+\phi)}|^2}\right]
+ 2 t\sin\theta\sin\phi \nonumber \\ \nonumber \\ 
&\le& -\frac{2t\sin\theta\sin\phi}{|1-te^{i(\theta+\phi)}|^2}
+ 2t \sin\theta\sin\phi \nonumber \\ \nonumber \\ 
&=& 2t\sin\theta\sin\phi\, 
\frac{-2t\cos(\theta+\phi)+t^2}
{|1-te^{i(\theta+\phi)}|^2} \nonumber \\ \nonumber \\
&\ll& 
t \sin\theta\sin\phi\, \frac{\max(t,t^2)}{|1-te^{i(\theta+\phi)}|^2}\,.
\label{**2}
\end{eqnarray}
    
If $t\le 1/2$, then
$$
|1-te^{i(\theta+\phi)}|^2 \gg 1\,,
$$
and we obtain
\begin{equation}
D \ll t^2 \sin\theta\sin\phi
\ll \eta t^2 + \eta^{-1} t^2 \sin^2\phi
=  \eta\, \frac{|z|^2}{|\zeta|^2}
+ \frac{|z|^2}{\eta} \left|{\rm Im}\frac{1}{\zeta} \right|^2\,,
\label{**3}
\end{equation}
with an arbitrary $\eta>0$.
  
If $t\ge 1/2$, then
$$ 
|1-te^{i(\theta+\phi)}|^2 \gg t^2 \sin^2\theta\,,
$$
so that (\ref{**2}) gives us
\begin{equation}
D \ll t \frac{\sin\phi}{\sin\theta}
\ll \frac{\eta}{\sin^2\theta} + \frac{t^2}{\eta} \sin^2\phi
=\frac{\eta}{\sin^2\theta} +
\frac{|z|^2}{\eta} \left|{\rm Im}\frac{1}{\zeta} \right|^2\,,
\label{**4}
\end{equation}
again, with an arbitrary positive $\eta$.
We shall use the bounds (\ref{**3}) and (\ref{**4}) when
$|z| \left|{\rm Im}\frac{1}{\zeta}  \right| \le 1$.

Now, for $z\in{\bf C}_+$, $r=|z|$, we have
\begin{eqnarray*}
u(z)-u_1(z) &=& \int_{\bf C_+} D(z,\zeta) \, d\mu(\zeta) \\ \\
&\le&
\left(\int_{\left|{\rm Im}\frac{1}{\zeta} \right|\ge \frac{1}{r}}   
+ \int_{\left|{\rm Im}\frac{1}{\zeta} \right|\le \frac{1}{r},
\,|\zeta|\ge 2r} +
\int_{\left|{\rm Im}\frac{1}{\zeta} \right|\le \frac{1}{r},
\, |\zeta|\le 2r}
\right) D(z,\zeta)\, d\mu(\zeta) \\ \\
&\ll& 
r \int_0^{r} \frac{d\n(t)}{t}
+ \frac{r^2}{\eta} \int_{r}^\infty \frac{d\n(t)}{t^2}
+ \frac{\eta}{\sin^2\theta} \int_0^{r} d\mu (t)
+ \eta r^2 \int_{r}^\infty \frac{d\mu (t)}{t^2} \\ \\
&\ll& \frac{\delta(r)}{\eta} +
\frac{\eta r^2}{\sin^2\theta} \int_{r}^\infty
\frac{M(t,u)}{t^3}\, dt\,.
\end{eqnarray*}

Then, using estimate (\ref{3.4}) for $u_1^+(z)$ in the upper
half-plane, we obtain estimate (\ref{5.a}) for $0<\theta<\pi$. The same
argument applies for the lower half-plane, and the proof of (\ref{5.a}) is 
complete.

\subsection{Integral inequality for $M(r,u)$.}

Here we prove the integral inequality 
\begin{equation}
M(r,u) \ll
\sqrt{\delta^*(r)\, r^2 \int_r^\infty \frac{M(t,u)}{t^3}dt}\,.
\label{3.12}
\end{equation}

First, we improve estimate (\ref{5.a}) near the real axis. 
Consider the function $u(z)$ in the angles
$\{|\arg z|\le \pi/6\}$ and $\{|\arg z -\pi|\le \pi/6\}$.
On the boundary of these angles, 
$$
u(re^{i\theta}) \ll \Phi (r)\,,
\qquad \theta = \pm\frac{\pi}{6}, \quad  \pi \pm \frac{\pi}{6}\,,
$$
where
$$
\Phi (r) 
= \eta^{-1} \delta^*(r) 
+ \eta r^2 \int_r^\infty \frac{M(t,u)}{t^3}\, dt\,.
$$
Applying Lemma~1 to $u(z)$ in these angles, we obtain
for $|\theta|\le \pi/8$ and $|\pi -\theta|\le \pi/8$,
$$
u(re^{i\theta}) \ll
r^{-3} \int_0^r \Phi (t) t^2 dt 
+ r^3 \int_r^\infty \frac{\Phi(t)}{t^4}dt 
\ll \Phi (r)\,.
$$
The second inequality follows since the function $\Phi (r)$ does not
decrease, and the function $r^{-2}\Phi (r)$ does not increase.

Thus, for $0<r<\infty$, 
$$
M(r,u) \ll \Phi (r)= \eta^{-1} \delta^*(r) 
+ \eta r^2 \int_r^\infty \frac{M(t,u)}{t^3}\, dt\,.
$$
Choosing  
$$
\eta = 
\sqrt{\delta^*(r)}\, : \, 
\sqrt{r^2 \int_r^\infty \frac{M(t,u)}{t^3}dt}\,,
$$ 
we obtain inequality (\ref{3.12}).

\subsection{Solution of the integral inequality (\ref{3.12}).}
We set
$$
M_1(r) = \int_r^\infty \frac{M(t,u)}{t^3}\,dt\,.
$$
Then
$$
M(r,u) = -\,r^3 M_1'(r)\,,
$$
and inequality (\ref{3.12}) takes the form 
$$
-\,M_1'(r)r^2 \ll \sqrt{\delta^*(r)\,M_1(r)}
$$
or
$$
-\frac{d\sqrt {M_1(r)}}{dr} 
\ll \frac{\sqrt{\delta^*(r)}}{r^2}
\,.
$$
Integrating this inequality from $\infty$ to $r$,
we obtain
$$
M_1(r) \ll 
\left[
\int_r^\infty \frac{\sqrt{\delta^*(t)}}{t^2}\,dt 
\right]^2\,.
$$
On the other hand, since $M(r,u)$ does not decrease,
$$
M_1(r) \ge M(r,u) \int_r^\infty \frac{dt}{t^3}
= \frac{M(r,u)}{2r^2}\,.
$$
Therefore, 
$$
M(r,u) \ll r^2 M_1(r) \ll r^2 \left[
\int_r^\infty \frac{\sqrt{\delta^*(t)}}{t^2}\,dt 
\right]^2\,,
$$
completing the proof of Theorem~2.
$\Box$

\section{Proof of Theorem 3}
\setcounter{equation}{0}

We divide the proof into 4 parts. 
Set 
$$
B:= \limsup_{r\to 0} \frac{\mu(r)}{r}\,,
$$
$$
C:= \int_0^\infty \frac{\delta (t)}{t^2}\, dt\,.
$$
Without loss of generality, we assume that both values $B$ and $C$ are
finite.

First, we shall prove the theorem under the additional assumption 
\begin{equation}
{\rm supp}(\mu) \subset \bf{\bar C}_-\,,
\label{**}
\end{equation}
and till Section~6.4 we assume that the function $u(z)$ is
harmonic in ${\bf C}_+$.

\subsection{The function $u(z)$ 
has nonnegative harmonic majorants in $\bf C_\pm$.}
Consider the function
$$
U(z):=-u(z)-\frac{y}{\pi}\int_{-\infty}^\infty 
\frac{u^-(t)dt}{(t-x)^2+y^2}\,.
$$
This function is harmonic in ${\bf C}_+$ and $U(x)\le 0,\; x\in{\bf R}$.
Moreover, for $y>0$,
\begin{eqnarray*}
U(iy)\le -u(iy)&=&-\lim\limits_{\varepsilon\to
0}\int_{|\zeta|\ge\varepsilon, \,
\zeta\in\bf{\bar C}_-}\left[\log
\left|1-\frac{iy}{\zeta}\right| 
+ \mbox{Re}\,\frac{iy}{\zeta}\right]d\mu(\zeta)
\\ \\ 
&\le& -\lim\limits_{\varepsilon\to 0}\int_{|\zeta|\ge\varepsilon,\,
\zeta\in
\bf{\bar C}_-}\mbox{Re}\,\frac{iy}{\zeta}\,d\mu(\zeta)
\\ \\
&=&y\int_{\bf{\bar C}_-}\mbox{Im}\,
\frac{1}{\zeta}\,d\mu(\zeta)\le Cy\,.
\end{eqnarray*}

By the Poisson-Nevanlinna representation of harmonic functions in the 
semi-disk $D_R$ (cf. Section~3), we have
\begin{eqnarray}
U(z)&\le&
\frac{1}{2\pi}\int_0^{\pi}U(Re^{i\phi}) K_2(z,Re^{i\phi})\, d\phi 
\nonumber \\ \nonumber \\
&\le&\frac{2Rr(R+r)}{\pi(R-r)^3}\int_0^{\pi}u^-(Re^{i\phi})\,d\phi
\nonumber \\ \nonumber \\
&\le&
\frac{4Rr(R+r)}{(R-r)^3}\,T(R,u)\,,
\label{1}
\end{eqnarray}
where $T(R,u)$
is the Nevanlinna characteristic of $u$.

Note that, for any $\delta>0$, the function $u$ can be
represented in the form
\begin{eqnarray}
u(z)&=&\int_{|\zeta|>\delta}\left[\log\left|1-\frac{z}{\zeta}\right|+\mbox{Re}
\frac{z}{\zeta}\right]\,d\mu(\zeta) \nonumber \\
\nonumber \\
&\ & \qquad +  \int_{|\zeta|\le \delta}\log\left|1-
\frac{z}{\zeta}\right|\,d\mu(\zeta)+\mbox{Re}\left(z
\int_{|\zeta|<\delta}\frac{d\mu(\zeta)}{\zeta}\right) \nonumber \\
\nonumber \\
&=:&u^{\delta}(z)+v^{\delta}(z)+\mbox{Re}\left(z
\int_{|\zeta|<\delta}\frac{d\mu(\zeta)}{\zeta}\right).
\label{2}
\end{eqnarray}
The well-known Borel estimates
$$\max_{|z|\le r}u^{\delta}(z)=o(|z|^2), \qquad \max_{|z|\le r}v^{\delta}(z)=o(|z|),\quad
z\to\infty,$$
imply that $T(R,u)=o(R^2),\; R\to\infty$. Therefore, 
by setting $R=2r$ in (\ref{1}),
we get
$$U^+(z)=o(|z|^2),\qquad z\to\infty,\; \qquad \mbox{Im}z>0.$$

Applying the Phragm\'{e}n-Lindel\"{o}f principle in the angles
$\{ 0<\arg z< \pi/2\}$ and $\{ \pi/2<\arg z< \pi \}$, we
conclude that
\begin{equation}
U(z)\le Cy,\qquad z=x+iy\in {\bf C_+}\,;
\label{3}
\end{equation}
i.e. $-U(z)+Cy$ is a nonnegative harmonic function in ${\bf C_+}$. 
Since $u(z) \le -U(z)+Cy$,
the function $u(z)$ also has a nonnegative harmonic majorant in ${\bf
C_+}$.

For $z\in{\bf C}_-$, we write
\begin{eqnarray}
u(z)-u(\bar{z}) &=&
\lim\limits_{\varepsilon\to 0}\,
\int_{|\zeta|\ge\varepsilon,\zeta\in{\bf C}_-}
\log\left|\frac{1-z/\zeta}{1-z/\bar{\zeta}}\right|+\mbox{Re}\left[z\left(
\frac{1}{\zeta}-\frac{1}{\bar{\zeta}}\right)\right]\,d\mu(\zeta)
\nonumber \\ \nonumber \\ 
&\le& 
\int_{{\bf C}_-}\mbox{Re}\left[z2i\mbox{Im}\frac{1}{\zeta}\right]\, 
d\mu(\zeta) \le 2C|y|\,.
\label{4}
\end{eqnarray}
Because $u(\bar{z})$ has a nonnegative harmonic majorant in ${\bf C}_-$, we
get the desired conclusion.

\subsection{Estimate of $u(z)$ near the origin.} Set
$$I(r):=\frac{1}{r}\int_0^{\pi}u(re^{i\theta})\sin\theta\,d\theta.$$
Let us prove that
\begin{equation}
\limsup\limits_{r\to 0}I(r)\ll B\,.
\label{5}
\end{equation}
 
For any given $\varepsilon>0$, choose a positive $\delta<\varepsilon$ such
that
$$\mu(r)<(B+\varepsilon)r,\quad\mbox{for}\quad 0<r<\delta.$$
Let us represent $u$ by the formula (\ref{2}) with this $\delta$. Since
$$u^{\delta}(z)=O(|z|^2),\quad z\to 0,$$
we have
\begin{equation}
\limsup\limits_{r\to 0}I(r)\le\limsup\limits_{r\to 0}I^{\delta}(r)+
\left|\int_{|\zeta|<\delta}\frac{d\mu(\zeta)}{\zeta}\right|\,,
\label{6}
\end{equation}
where
$$I^{\delta}(r)=\frac{1}{r}\int_0^{\pi}v^{\delta}(re^{i\theta})\sin\theta\,
d\theta.$$
It suffices to show that
\begin{equation}
I^{\delta}(r)\ll  B+\varepsilon+
\int_{|\zeta|<\delta}\mbox{Im}\frac{1}{\zeta}\,d\mu(\zeta)\,, \qquad
0<r<\delta. 
\label{7}
\end{equation}
Indeed, if (\ref{7}) is valid, then substituting it into (\ref{6}), we get
$$
\limsup\limits_{r\to 0}I(r) \ll B+\varepsilon+
\int_{|\zeta|<\delta}\mbox{Im}\frac{1}{\zeta}\,d\mu(\zeta)
+\left|\int_{|\zeta|<\delta}\frac{d\mu(\zeta)}{\zeta}\right|.$$
Taking the limit as $\varepsilon\to 0$ (then $\delta\to 0$ as well), we
obtain
(\ref{5}).

To prove (\ref{7}), we set for $|z|=r,\; 0<r<\delta$:
$$
v^{\delta}(z)=\int_{r<|\zeta|<\delta}+\int_{|\zeta|<r}=:v_1^{\delta}(z)
+v_2^{\delta}(z)\,,
$$
and
$$
I_j^{\delta}(r):=\frac{1}{r}\int_0^{\pi}v_j^{\delta}(re^{i\theta})
\sin\theta\,d\theta,\qquad j=1,2\,.
$$
Note that
\begin{eqnarray*}
\int_0^{\pi}\left|\log\left|1-\frac{re^{i\theta}}{\zeta}\right|\right|
\sin\theta\,d\theta
&\le&
2\int_0^{2\pi}\log^+\left|1-\frac{re^{i\theta}}{\zeta}\right|\,d\theta
-\int_0^{2\pi}\log\left|1-\frac{re^{i\theta}}{\zeta}\right|\,d\theta \\ \\ 
&\le&
4\pi\log\left(1+\frac{\delta}{|\zeta|}\right)+2\pi\log^+\frac{r}{|\zeta|}\,.
\end{eqnarray*}
This estimate will allow us to change the integration order in the double
integrals that arise when estimating $I_j^{\delta}(r),\; j=1,2,$ below.

Write
$$
I_1^{\delta}(r)=\frac{1}{r}\int_{r<|\zeta|<\delta}d\mu(\zeta)\int_0^{\pi}\log
\left|1-\frac{re^{i\theta}}{\zeta}\right|\sin\theta\,d\theta\,.
$$
For $r<|\zeta|$, we have
\begin{eqnarray*}
\int_0^{\pi}\log\left|1-\frac{re^{i\theta}}{\zeta}\right|\sin\theta\,d\theta
&=&
-\mbox{Re}\sum_{k=1}^{\infty}\frac{r^k}{k\zeta^k}\int_0^{\pi}
e^{ik\theta}\sin\theta\,d\theta \\ \\ 
&=&
\frac{\pi r}{2}\mbox{Im}\frac{1}{\zeta}+\mbox{Re}\sum_{m=1}^{\infty}
\frac{r^{2m}}{m(4m^2-1)\zeta^{2m}}\,.
\end{eqnarray*}
Hence
$$
I_1^{\delta}(r)\le\frac{\pi}{2}\int_{|\zeta|<\delta}\mbox{Im}\frac{1}{\zeta}\,
d\mu(\zeta)+
\sum_{m=1}^{\infty}\frac{r^{2m-1}}{m(4m^2-1)}\int_{r<|\zeta|<\delta}
\frac{d\mu(\zeta)}{|\zeta|^{2m}}\,.
$$
Since
\begin{eqnarray*}
\int_{r<|\zeta|<\delta}\frac{d\mu(\zeta)}{|\zeta|^{2m}}&=&
\int_r^{\delta}\frac{d\mu(t)}{t^{2m}} \\ \\
&\le&
\frac{\mu(\delta)}{\delta^{2m}}
+2m\int_r^{\delta}\frac{\mu(t)}{t^{2m+1}}\,dt \\ \\
&\le& (B+\varepsilon)\delta^{-2m+1} + 2(B+\varepsilon)r^{-2m+1}
< 3(B+\varepsilon)r^{-2m+1}\,,
\end{eqnarray*}
we get
\begin{equation}
I_1^{\delta}(r)\le\frac{\pi}{2}\int_{|\zeta|<\delta}\mbox{Im}
\frac{1}{\zeta}\,d\mu(\zeta) + 3(B+\varepsilon)
\sum_{m=1}^{\infty}\frac{1}{m(4m^2-1)}\,.
\label{8}
\end{equation}

Further,
$$I_2^{\delta}(r)=\frac{1}{r}\int_{|\zeta|<r}d\mu(\zeta)\int_0^{\pi}\log\left|
1-\frac{re^{i\theta}}{\zeta}\right|\sin\theta\,d\theta.$$
For $|\zeta|\le r$, we have
\begin{eqnarray*}
\int_0^{\pi}\log\left|1-\frac{re^{i\theta}}{\zeta}\right|
\sin\theta\,d\theta
&=&
2\log\frac{r}{|\zeta|}+\int_0^{\pi}
\log\left|1-\frac{\bar{\zeta}e^{i\theta}}{r}\right|\sin\theta\,d\theta
\\ \\ 
&=&
2\log\frac{r}{|\zeta|}-\mbox{Re}\sum_{k=1}^{\infty}\frac{\bar{\zeta}^k}
{kr^k}\int_0^{\pi}e^{ik\theta}\sin\theta\,d\theta \\ \\
&=& 2\log\frac{r}{|\zeta|}-\frac{\pi}{2r}\mbox{Im}\zeta+\mbox{Re}
\sum_{m=1}^{\infty}\frac{\bar{\zeta}^{2m}}{m(4m^2-1)r^{2m}}\,.
\end{eqnarray*}
Whence
\begin{eqnarray}
I_2^{\delta}(r) &\le&
\frac{1}{r}\int_{|\zeta|<r}\left[2\log\frac{r}{|\zeta|}
+\frac{\pi}{2}+
\sum_{m=1}^{\infty}\frac{1}{m(4m^2-1)}\right]\,d\mu(\zeta)
\nonumber \\ \nonumber \\
&\ll&
\frac{1}{r}\int_0^r \log\frac{r}{t}d\mu(t)+ \frac{\mu(r)}{r} 
\nonumber \\ \nonumber \\
&\ll& B+\varepsilon 
\label{9}
\end{eqnarray}
Since $I^{\delta}=I_1^{\delta}+I_2^{\delta}$, the desired inequality (\ref{7})
follows from (\ref{8}) and (\ref{9}).

\subsection{Estimate of $u^+(z)$ on the real and imaginary axes.} Let us
prove that
\begin{equation}
\int_{-\infty}^\infty 
\frac{u^+(t)}{t^2}dt+\limsup\limits_{y\to+\infty}\frac
{u^+(iy)}{y}\ll B+C\,.
\label{10'} 
\end{equation}
Since $u$ has a nonnegative harmonic majorant in ${\bf C}_+$, we have
$$
\int_{-\infty}^{\infty}\frac{|u(t)|dt}{1+t^2}<\infty\,,
$$
and $u$ admits the Poisson representation
\begin{equation}
u(re^{i\varphi})=\frac{r\sin\varphi}{\pi}
\int_{-\infty}^{\infty}\frac{u(t)\,dt}
{r^2+t^2-2rt\cos\varphi}+kr\sin\varphi,
\label{11}
\end{equation}
where
$$k=\limsup\limits_{y\to+\infty}\frac{u(iy)}{y}\neq\infty.$$
Note that inequality (\ref{3}) implies
$$
u(iy)\ge-\frac{y}{\pi}\int_{-\infty}^{\infty}\frac{u^-(t)dt}{t^2+y^2}-Cy=
o(y)-Cy, \qquad y\to\infty\,;
$$
i.e. $k\ge -C$.

Multiplying (\ref{11}) by $\sin\varphi$, integrating against $\varphi$
from $0$ to $\pi$, and taking into account that
$$\int_0^{\pi}\frac{\sin^2\varphi\,d\varphi}{r^2+t^2-2rt\cos\varphi}=
\frac{\pi}{2}\min\left(\frac{1}{r^2},\frac{1}{t^2}\right),$$
we get
$$\int_0^{\pi}u(re^{i\varphi})\sin\varphi\,d\varphi=
\frac{r}{2}\int_{-\infty}^{\infty}u(t)\min\left(\frac{1}{r^2},\frac{1}{t^2}
\right)\,dt+\frac{k\pi r}{2}.$$
Hence
$$\int_{-\infty}^{\infty}u^+(t)\min\left(\frac{1}{r^2},\frac{1}{t^2}\right)
\,dt+k^+\pi
=\int_{-\infty}^{\infty}u^-(t)\min\left(\frac{1}{r^2},
\frac{1}{t^2}\right)\,dt+k^-\pi+2I(r).$$
Letting $r\to 0$, we obtain by the monotone convergence theorem
$$\int_{-\infty}^{\infty}\frac{u^+(t)}{t^2}\,dt+k^+\pi=
\int_{-\infty}^{\infty}\frac{u^-(t)}{t^2}\,dt+k^-\pi+2\lim\limits_{r\to 0}
I(r)$$
(it turns out that the last limit exists). Taking into account that
$k^-\le C$ and using (\ref{5}), we get (\ref{10'}).

\subsection{Concluding steps.} From (\ref{3}) and (\ref{10'}) we obtain
$$
\limsup\limits_{y\to-\infty}\frac{u^+(iy)}{|y|}\le\limsup\limits_
{y\to+\infty}\frac{u^+(iy)+2Cy}{y}\ll B+C\,.
$$
Since $u$ has nonnegative harmonic majorants in both upper and lower
half-planes, the following inequality holds in the whole plane:
\begin{equation}
u(z)\ll \frac{|y|}{\pi}\int_{-\infty}^\infty 
\frac{u^+(t)\,dt}{(t-x)^2+y^2}
 + (B+C)|y|, \qquad z=x+iy.
\label{12}
\end{equation}
The assertion of Theorem~3 can be obtained from this inequality and 
(\ref{10'}) by applying a known argument (cf. \cite{MS1}, \cite{MS3}). 
First, one applies (\ref{12}) and (\ref{10'}) to get the upper bound for 
$u(z)$ in the angles $\{|\arg z \pm \pi/2|\le \pi/4\}$, and then, 
using the Phragm\'en-Lindel\"of principle, one obtains the upper bound
for $u(z)$ in the complementary angles. 
This gives $u(z)\ll (B+C)|z|$, and completes the proof of estimate
(\ref{kolmogorov}) for the special case (\ref{**}).

Now, let $\mu$ be an arbitrary measure in $\bf C$ satisfying conditions of
Theorem~3 and having finite value $C$. As in Section~5, we define the
measure $\mu_1$, ${\rm supp}\mu_1\subset {\bf \bar C}_-$, 
by reflecting with respect to the real axis the part
of $\mu$ which charges $\bf C_+$.
Since
$$
\int_{\bf{C}_-} {\rm Im} \frac{1}{\zeta}\, d\mu_1 (\zeta) \le C <\infty\,,
$$
the measure $\mu_1$ also satisfies the conditions of Theorem~3, and we can
define the corresponding generalized canonical integral $u_1(z)$ of this
measure. Then a straightforward estimate (cf. \ref{4}) shows that for
$z\in \bf{\bar C}_+$
$$
u(z) \le u_1(z) + 2|y|\int_{\bf C} 
\left| {\rm Im} \frac{1}{\zeta} \right|\, d\mu(\zeta)
\ll (B+C)|z|\,.
$$
The same estimate holds in the lower half-plane, and the general case of
Theorem~3 follows. 
$\Box$

\bigskip\par\noindent Vladimir Matsaev: {\em School of Mathematical
Sciences,
Tel-Aviv University, \newline
\noindent Ramat-Aviv, 69978, Israel

\par\noindent matsaev@math.tau.ac.il} 

\bigskip\par\noindent Iossif Ostrovskii: {\em 
Department of Mathematics, Bilkent University, \newline
06533 Bilkent, Ankara, Turkey 

\par\noindent iossif@fen.bilkent.edu.tr}

\medskip\par\noindent and {\em Verkin Institute for Low Temperature
Physics and Engineering, \newline
310164 Kharkov, Ukraine

\par\noindent ostrovskii@ilt.kharkov.ua}

\bigskip\par\noindent Mikhail Sodin: {\em School of Mathematical
Sciences, Tel-Aviv University, \newline
\noindent Ramat-Aviv, 69978, Israel

\par\noindent sodin@math.tau.ac.il}

\end{document}